\documentclass[11pt]{article}

\usepackage{amssymb, amsmath, amsthm, graphicx}
\usepackage[left=1in,top=1in,right=1in]{geometry}

\usepackage{subfigure}

      \theoremstyle{plain}
      \newtheorem{theorem}{Theorem}[section]
      \newtheorem{lemma}[theorem]{Lemma}
            
            \newtheorem{observation}[theorem]{Observation}

      \theoremstyle{definition}

      \theoremstyle{remark}

\title{$k$-quasi planar graphs}

\author{Andrew Suk\thanks{EPFL, Lausanne. Email: {\tt suk@cims.nyu.edu}.  The author gratefully acknowledges the support from the Swiss National
Science Foundation, Grant No. 200021-125287/1.}}
\begin{document}

\maketitle

\begin{abstract}

A topological graph is \emph{$k$-quasi-planar} if it does not contain $k$ pairwise crossing edges.  A topological graph is \emph{simple} if every pair of its edges intersect at most once (either at a vertex or at their intersection).  In 1996, Pach, Shahrokhi, and Szegedy \cite{pach} showed that every $n$-vertex simple $k$-quasi-planar graph contains at most $O\left(n(\log n)^{2k-4}\right)$ edges.  This upper bound was recently improved (for large $k$) by Fox and Pach \cite{fox} to $n(\log n)^{O(\log k)}$.  In this note, we show that all such graphs contain at most $(n\log^2n )2^{\alpha^{c_k}(n)}$ edges, where $\alpha(n)$ denotes the inverse Ackermann function and $c_k$ is a constant that depends only on $k$.

\end{abstract}

\section{Introduction}

A \emph{topological graph} is a graph drawn in the plane such that its vertices are represented by points
 and its edges are represented by non-self-intersecting arcs connecting the corresponding points. The arcs are
allowed to intersect, but they may not pass through vertices except for
their endpoints.  Furthermore, the edges are not allowed to have tangencies, i.e., if two edges share an interior point, then they must properly cross at that point in common.
We only consider graphs without parallel edges or self-loops.
A topological graph is \emph{simple} if every pair of its edges intersect at most once.  If the edges are drawn as straight-line segments, then the graph is \emph{geometric}.  Two edges of a topological graph \emph{cross} if their interiors share a point.

Finding the maximum number of edges in a topological graph with a forbidden crossing pattern has been a classic problem in extremal topological graph theory (see \cite{ackerman2, acktardos, agarwal, cap, fox, fulek, pinchasi, tardos, valtr}).  It follows from Euler's Polyhedral Formula that every topological graph on $n$ vertices and no crossing edges has at most $3n-6$ edges.  A topological graph is $k$\emph{-quasi-planar}, if it does not contain $k$ pairwise crossing edges.  Hence 2-quasi-planar graphs are planar.  An old conjecture (see Problem 1 in section 9.6 of \cite{brass}) states that for any fixed $k > 0$, every $k$-quasi-planar graph on $n$ vertices has at most $c_kn$ edges, where $c_k$ is a constant that depends only on $k$.  Agarwal et al. were the first to prove this conjecture for simple 3-quasi-planar graphs.  Later Pach, Radoi\v{c}i\'c, and T\'oth \cite{rados} generalized the result for all (not simple) 3-quasi-planar graphs.  Recently, Ackerman \cite{ackerman} proved the conjecture for $k = 4$.

For $k \geq 5$, Pach, Shahrokhi, and Szegedy \cite{pach} showed that every simple $k$-quasi-planar graph on $n$ vertices has at most $c_kn(\log n)^{2k-4}$ edges.  This bound can be improved to $c_kn(\log n)^{2k-8}$ by using a result of Ackerman \cite{ackerman}.  Valtr \cite{valtrpar} proved that every $n$-vertex $k$-quasi-planar geometric graph contains at most $O(n\log n)$ edges.  Later, he extended this result to simple topological graphs with edges drawn as $x$-monotone curves \cite{valtr}.  Pach, Radoi\v{c}i\'c, and T\'oth showed that every $n$-vertex (not simple) $k$-quasi-planar graph has at most $c_kn(\log n)^{4k-12}$ edges, which can also be improved to $c_kn(\log n)^{4k-16}$ by a result of Ackerman \cite{ackerman}.

Very recently, Fox and Pach \cite{fox} improved (for large $k$) the exponent in the polylogarithmic factor for simple topological graphs.  They showed that every simple $k$-quasi-planar graph on $n$ vertices has at most $n(c\log n/\log k)^{c\log k}$ edges, where $c$ is an absolute constant.  Our main result is the following.

\begin{theorem}
\label{simple}
Let $G = (V,E)$ be an $n$-vertex simple $k$-quasi-planar graph.  Then $|E(G)| \leq   (n\log^2 n)2^{\alpha^{c_k}(n)}$, where $\alpha(n)$ denotes the inverse Ackermann function and $c_k$ is a constant that depends only on $k$.
\end{theorem}

In the proof of Theorem \ref{simple}, we apply results on generalized Davenport-Schinzel sequences.  This method was used by Valtr \cite{valtr}, who showed that every $n$-vertex simple $k$-quasi-planar graph with edges drawn as $x$-monotone curves has at most $2^{2^{ck}}n\log n$ edges, where $c$ is an absolute constant.  Our next theorem extends his result to (not simple) topological graphs with edges drawn with $x$-monotone curves, and moreover we obtain a slightly better upper bound.

\begin{theorem}
\label{xmono}
Let $G =(V,E)$ be an $n$-vertex (not simple) $k$-quasi planar graph with edges drawn as $x$-monotone curves.  Then $|E(G)| \leq 2^{ck^{3}}n\log n$, where $c$ is an absolute constant.

\end{theorem}

\section{Generalized Davenport-Schinzel sequences}

The sequence $u = a_1,a_2,...,a_m$ is called $l$-regular if any $l$ consecutive terms are pairwise different.  For integers $l,t \geq 2$, the sequence

$$S = s_1,s_2,...,s_{l t}$$

\noindent of length $l\cdot t$ is said to be of type $up(l,t)$ if the first $l$ terms are pairwise different and for $i = 1,2,...,l$

$$s_i = s_{i +l } = s_{i + 2l} = \cdots = s_{i + (t-1)l }.$$

\noindent For example,

$$a,b,c,a,b,c,a,b,c,a,b,c,$$

\noindent would be an $up(3,4)$ sequence.  By applying a theorem of Klazar on generalized Davenport-Schinzel sequences, we have the following.

\begin{theorem}[\cite{klazar92}]
\label{klazar}
For $l\geq 2$ and $t \geq 3$, the length of any $l$-regular sequence over an $n$-element alphabet that does not contain a subsequence of type $up(l,t)$ has length at most

$$n\cdot l 2^{(lt-3)} \cdot (10l)^{10\alpha^{lt}(n)}.$$

\end{theorem}

\noindent  For $l \geq 2$, the sequence

$$S = s_1,s_2,...,s_{3l-2}$$

\noindent of length $3l-2$ is said to be of type \emph{up-down-up}$(l)$, if the first $l$ terms are pairwise different, and for $i = 1,2,...,l$,

$$s_i = s_{2l-i} = s_{(2l-2)+ i}.$$

\noindent For example,

$$a,b,c,d,c,b,a,b,c,d,$$

\noindent would be an \emph{up-down-up}$(4)$ sequence.  Valtr and Klazar showed the following.

\begin{lemma}[\cite{val}]
\label{updownup}
For $l\geq 2$, the length of any $l$-regular sequence over an $n$-element alphabet containing no subsequence of type up-down-up$(l)$ has length at most $2^{O(l)}n$.

\end{lemma}

\noindent For more results on generalized Davenport-Schinzel sequences, see \cite{gabriel, seth, pettie}.

\section{Simple topological graphs}

In this section, we will prove Theorem \ref{simple}.  For any partition of $V(G)$ into two disjoint parts, $V_1$ and $V_2$, let $E(V_1,V_2)$ denote the set of edges with one endpoint in $V_1$ and the other endpoint in $V_2$.   The {\it bisection width} of a graph $G$, denoted by $b(G)$, is the smallest nonnegative integer such that there is a partition of  the vertex set $V=V_1 \, \dot{\cup} \, V_2$ with $\frac{1}{3}\cdot |V|\leq |V_i|\leq \frac{2}{3}\cdot |V|$ for $i=1,2$, and  $|E(V_1,V_2)|= b(G)$.  We will use the following result by Pach et al.

\begin{lemma}[\cite{pach}]
\label{bisection}
If $G$ is a graph with $n$ vertices of degrees $d_1,...,d_n$, then

$$b(G) \leq 7cr(G)^{1/2}  + 2\sqrt{\sum\limits_{i = 1}^n d_i^2},$$

\noindent where $cr(G)$ denotes the crossing number of $G$.

\end{lemma}

\noindent Since $\sum_{i = 1}^n d^2_i \leq 2n|E(G)|$ holds for every graph, we have

\begin{equation}
\label{bisect}
b(G) \leq 7cr(G)^{1/2} + 3\sqrt{|E(G)|n}.
\end{equation}

\medskip

\noindent \textbf{Proof of Theorem \ref{simple}.}  Let $k \geq 5$ and $f_k(n)$ denote the maximum number of edges in a simple $k$-quasi-planar graph on $n$ vertices.  We will prove that

$$f_k(n) \leq  (n\log^2 n)2^{\alpha^{c_k}(n)}$$

\noindent where $c_k=  10^5\cdot 2^{k^2 + 2k}$.  For sake of clarity, we do not make any attempts to optimize the value of $c_k$.  We proceed by induction on $n$.  The base case $n < 7$ is trivial.  For the inductive step $n \geq 7$, let $G = (V,E)$ be a simple $k$-quasi-planar graph with $n$ vertices and $m = f_k(n)$ edges, such that the vertices of $G$ are labeled 1 to $n$.  The proof splits into two cases.

\medskip

\noindent \emph{Case 1.}  Suppose that $cr(G) \leq m^2/(10^4\log^2n)$.  By (1), there is a partition $V(G) = V_1\cup V_2$ with $|V_1|,|V_2| \leq 2n/3$ and the number of edges with one vertex in $V_1$ and one vertex in $V_2$ is at most

$$b(G) \leq 7cr(G)^{1/2} + 3\sqrt{mn} \leq 7\frac{m}{100\log n} + 3\sqrt{mn}.$$

\noindent Let $n_1 = |V_1|$ and $n_2 = |V_2|$.  Now if $7m/(100\log n) \leq 3\sqrt{mn}$, then we have

$$m \leq 43n\log^2n$$

\noindent and we are done since $\alpha(n) \geq 2$ and $k \geq 5$.  Therefore, we can assume $7m/(100\log n) > 3\sqrt{mn}$, which implies

\begin{equation}
b(G) \leq \frac{m}{7\log n}.
\end{equation}

By the induction hypothesis and equation (2), we have

$$
  \begin{array}{ccl}
   m & \leq & f_k(n_1) + f_k(n_2) + b(G)\\\\
  & \leq & \left( n_1\log^2(2n/3)\right)2^{\alpha^{c_k}(n)} + \left(n_2\log^2(2n/3)\right)2^{\alpha^{c_k}(n)} + b(G) \\\\
     & \leq & \left(n\log^2(2n/3)\right) 2^{\alpha^{c_k}(n)} +   \frac{m}{7\log n}  \\\\
     & \leq & (n\log^2 n)2^{\alpha^{c_k}(n)} - 2n2^{\alpha^{c_k}(n)}\log n\log(3/2) + n2^{\alpha^{c_k}(n)}\log^2(3/2) + \frac{m}{7\log n} \\\\
  \end{array}
 $$

\noindent which implies

$$m\left( 1 - \frac{1}{7\log n} \right) \leq (n\log^2n)2^{\alpha^{c_k}(n)}\left(1 - \frac{2\log(3/2)}{\log n} + \frac{\log^2(3/2)}{\log^2n}\right).$$

\noindent Hence

$$m\leq (n\log^2n)2^{\alpha^{c_k}(n)}\frac{1 - 2\log(3/2)\log^{-1} n +  \log^2(3/2)\log^{-2}n }{1 - 1/(7\log n)} \leq (n\log^2n)2^{\alpha^{c_k}(n)}.$$

\medskip

\noindent \emph{Case 2.}  Now suppose that $cr(G) \geq m^2/(10^4\log^2 n)$.  By a simple averaging argument, there exists an edge $e = uv$ such that at least $2m/(10^4\log^2n)$ other edges cross $e$.  Fix such an edge $e = uv$, and let $E'$ denote the set of edges that cross $e$.

\begin{figure}[h]
\begin{center}
\includegraphics[width=200pt]{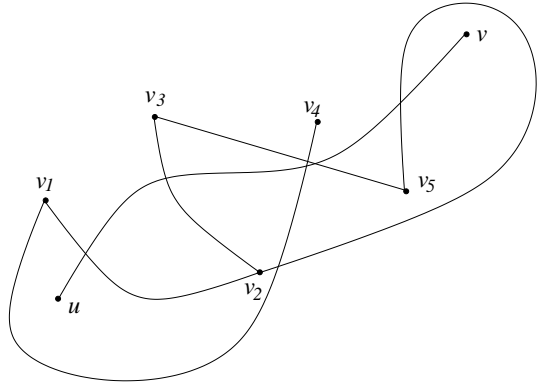}
  \caption{In this example, $S_1 = v_1,v_3,v_4,v_3,v_2$ and $S_2 = v_2,v_2,v_1,v_5,v_5$.}
  \label{s1s2}
 \end{center}
\end{figure}

We order the edges in $E'= \{e_1,e_2,...,e_{|E'|}\}$, in the order that they cross $e$ from $u$ to $v$.  Now we create two sequences $S_1 = a_1,a_2,...,a_{|E'|}$ and $S_2 = b_1,b_2,...,b_{|E'|}$ as follows.  For each $e_i \in E'$, as we move along edge $e$ from $u$ to $v$ and arrive at edge $e_i$, we turn left and move along edge $e_i$ until we reach its endpoint $u_i$.  Then we set $a_i = u_i$.  Likewise, as we move along edge $e$ from $u$ to $v$ and arrive at edge $e_i$, we turn right and move along edge $e_i$ until we reach its other endpoint $v_i$.  Then set $b_i = v_i$.  Thus $S_1$ and $S_2$ are sequences of length $|E'|$ over the alphabet $\{1,2,...,n\}$.  See Figure \ref{s1s2} for a small example.

Now we need the following two lemmas.  The first one is due to Valtr.

\begin{lemma}[\cite{valtr}]
\label{regular}
For $l \geq 1$, at least one of the sequences $S_1,S_2$ defined above contains an $l$-regular subsequence of length at least $|E'|/(4l)$.
\end{lemma}

$\hfill\square$

\begin{lemma}
\label{key}
Neither of the sequences $S_1$ nor $S_2$ contains a subsequence of type $up(2^{k^2 + k}, 2^k)$.

\end{lemma}

\noindent \textbf{Proof.}  By symmetry, it suffices to show that $S_1$ does not contain a subsequence of type $up(2^{k^2 + k},2^k)$.  We will prove by induction on $k$, that such a sequence will produce $k$ pairwise crossing edges in $G$.  The base cases $k = 1,2$ are trivial.  Now assume the statement holds up to $k-1$.  Let

$$S = s_1,s_2,...,s_{2^{k^2 + 2k}}$$

\noindent be our $up(2^{k^2 + k}, 2^k)$ sequence of length $2^{k^2 + 2k}$ such that the first $2^{k^2 + k}$ terms are pairwise different, and for $i = 1,2,...,2^{k^2 + k}$

$$s_i = s_{i + 2^{k^2 + k} }  = s_{i + 2\cdot 2^{k^2 + k} } = s_{i + 3\cdot 2^{k^2 + k}}= \cdots = s_{i + (2^k - 1)2^{k^2 + k}}.$$

\noindent For each $i = 1,2,...,2^{k^2 + k}$, let $v_i\in V_1$ denote the label (vertex) of $s_i$.  Moreover, let $a_{i,j}$ be the arc emanating from vertex $v_i$ to the edge $e$ corresponding to $s_{i + j2^{k^2 + k}}$ for $j = 0,1,2,...,2^k - 1$.  We will think of $s_{i + j2^{k^2 + k}}$ as a point on $a_{i,j}$ very close but not on edge $e$.  For simplicity, we will let $s_{2^{k^2 + 2k} + t} = s_t$ for all $ t\in \mathbb{N}$ and $a_{i,j} = a_{i,j\mod 2^k}$ for all $j  \in \mathbb{Z}$.  Hence there are $2^{k^2 + k}$ distinct vertices $v_1,...,v_{2^{k^2 + k}}$, each vertex of which has $2^k$ arcs emanating from it to the edge $e$.

\begin{figure}[h]
\begin{center}
\includegraphics[width=200pt]{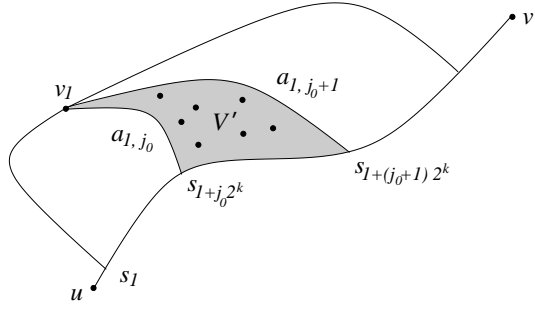}
  \caption{Vertices of $V'$ lie in the region enclosed by $a_{1,j_0},a_{1,j_0 + 1},e$.}
  \label{between}
 \end{center}
\end{figure}

Consider the drawing of the $2^k$ arcs emanating from $v_1$ and the edge $e$.  This drawing partitions the plane into $2^k$ regions.  By the Pigeonhole principle, there is a subset $V' \subset \{v_1,...,v_{2^{k^2 + k}}\}$ of size

$$\frac{2^{k^2 + k} - 1}{2^k},$$

\noindent such that all of the vertices of $V'$ lie in the same region.  Let $j_0 \in\{0,1,2,...,2^k-1\}$ be an integer such that $V'$ lies in the region bounded by $a_{1,j_0},a_{1,j_0 + 1}, e$.  See Figure \ref{between}.  In the case $j_0 = 2^k - 1$, $V'$ lies in the unbounded region.

 Let $v_i \in V'$ and $a_{i,j_0 + j_1}$ be an arc emanating out of $v_i$ for $j_1 \geq 1$.  Notice that $a_{i,j_0 + j_1}$ cannot cross both $a_{1,j_0}$ and $a_{1,j_0+1}$ since $G$ is simple.  Suppose that $a_{i,j_0 + j_1}$ crosses $a_{1,j_0 + 1}$.  Then the set of arcs (emanating out of $v_i$)

 $$A = \{a_{i,j_0 + 1}, a_{i,j_0 + 2},...,a_{i,j_0 + j_1-1}\}$$

 \noindent must also cross $a_{1,j_0 + 1}$.  Indeed, let $\gamma$ be the simple closed curve created by the arrangement

 $$a_{i,j_0 + j_1}\cup a_{1,j_0 + 1}\cup e.$$

 \noindent Since $a_{i,j_0 + j_1},a_{1,j_0 + 1},e$ pairwise intersect at precisely one point, $\gamma$ is well defined.  We define points $x = a_{i,j_0 + j_1} \cap a_{1,j_0  + 1}$ and $y = a_{1,j_0 + 1}\cap e$, and orient $\gamma$ in the direction from $x$ to $y$ along $\gamma$.

\begin{figure}[h]
    \label{gamma}
  \centering
  \subfigure[The case when $ j_0 + j_1 \mod 2^k \leq 2^{k} - 1$. ]{\label{j1p}\includegraphics[width=0.45\textwidth]{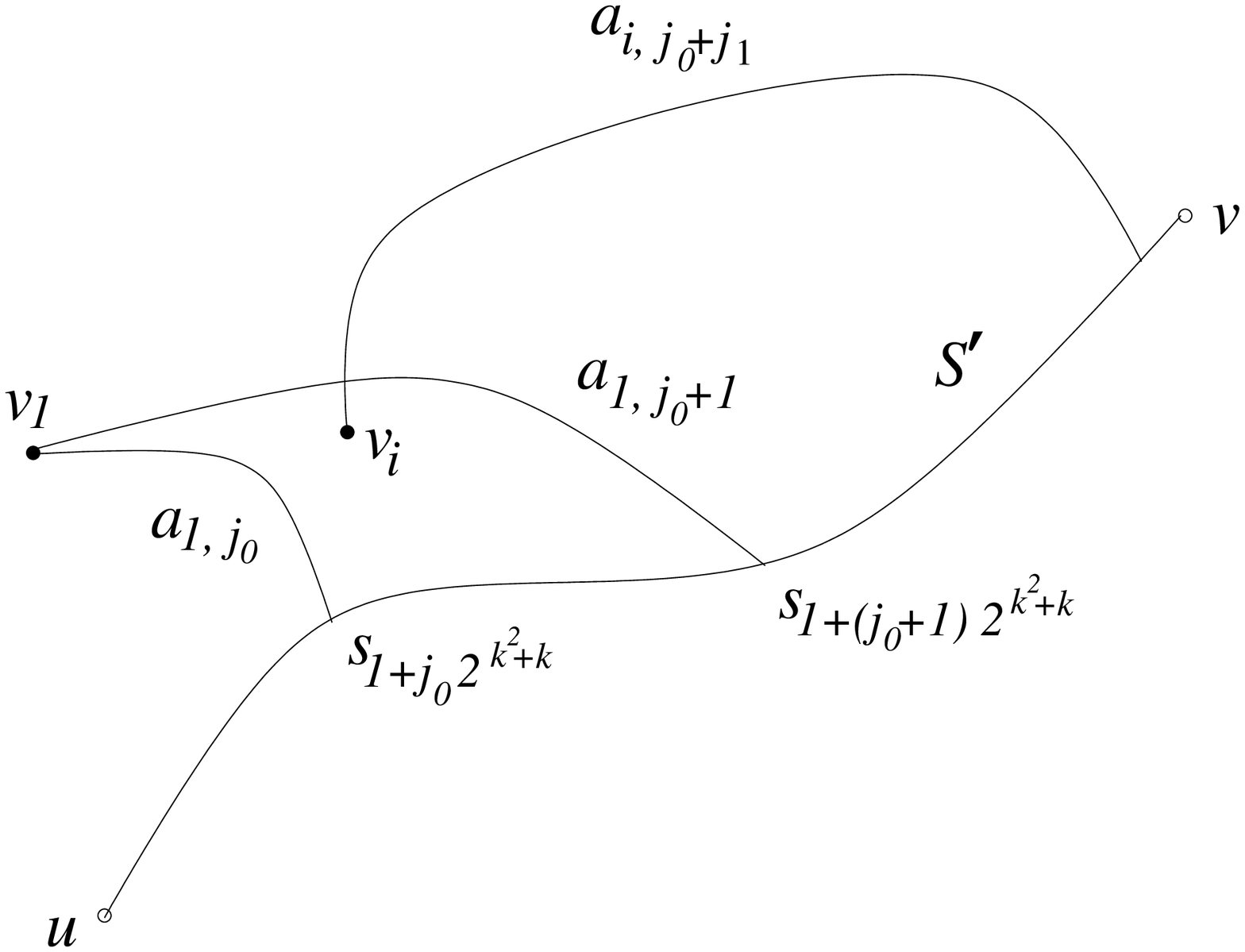}} \hspace{1cm}
\subfigure[$\gamma$ defined from Figure \ref{j1p}.]{\label{j1}\includegraphics[width=0.45\textwidth]{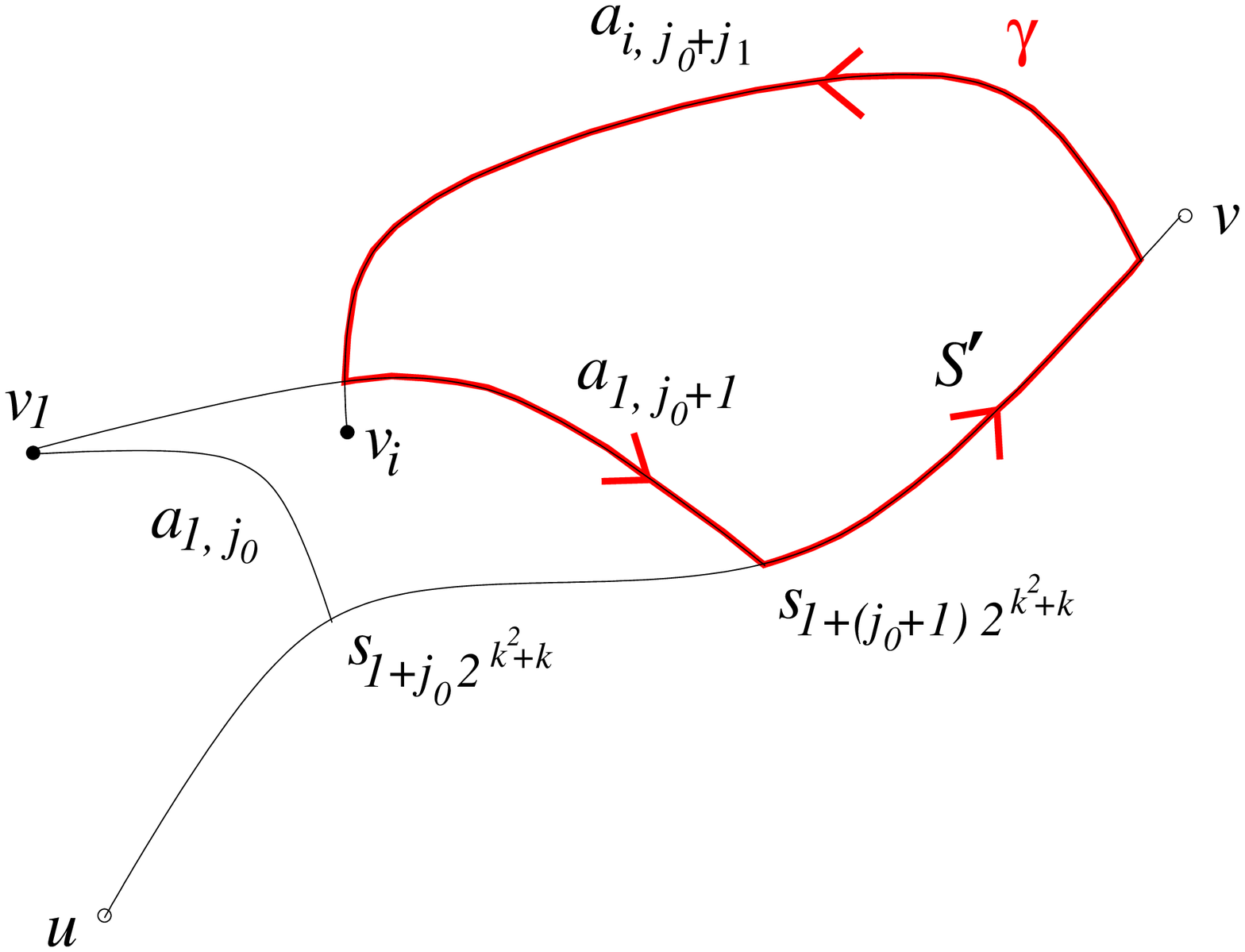}}  \hspace{1cm}
  \subfigure[The case when $ j_0 + j_1 \mod 2^k  < j_0$.  Recall $a_{i,j_0 + j_1} = a_{i,j_0+ j_1\mod 2^k}$.]{\label{j13p}\includegraphics[width=0.45\textwidth]{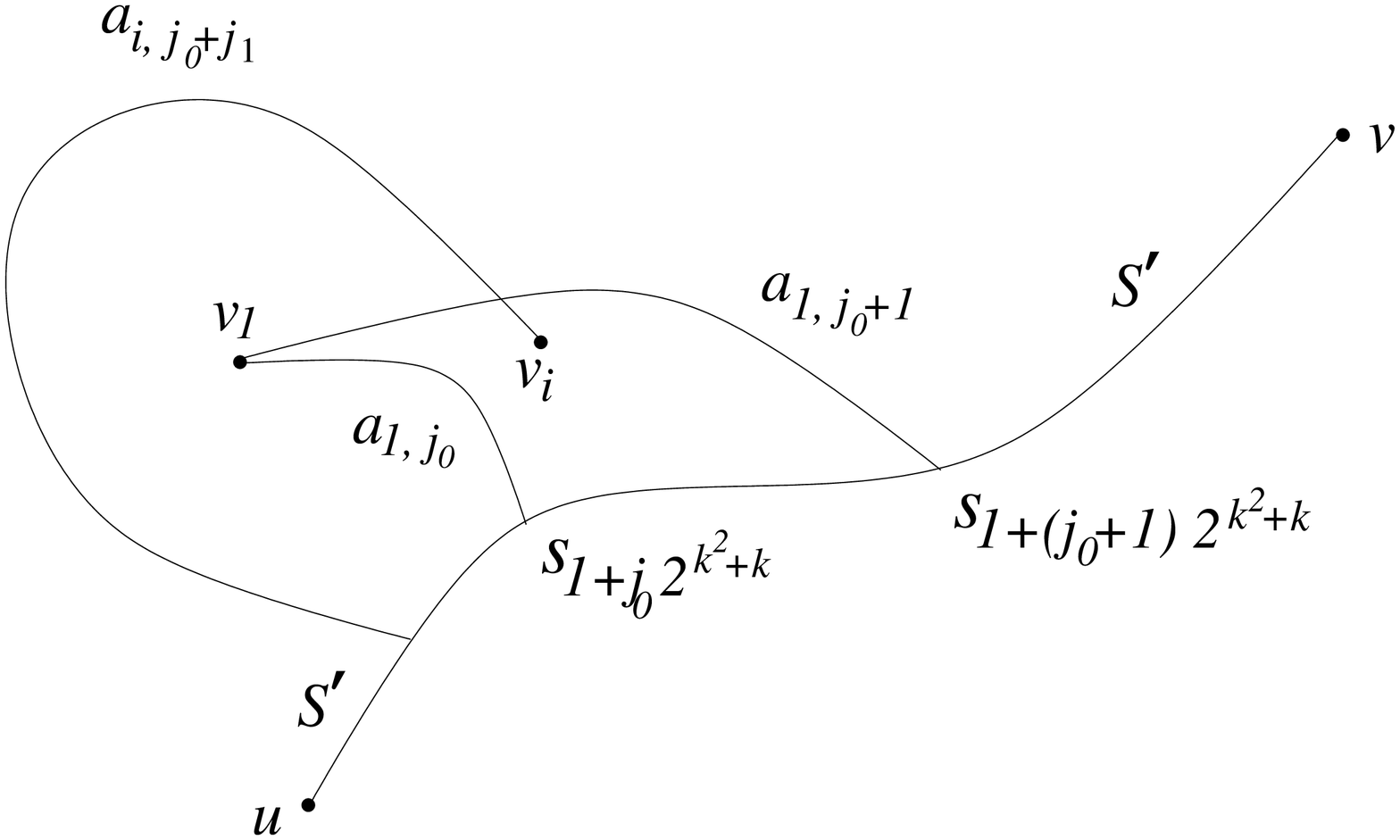}} \hspace{1cm}
\subfigure[$\gamma$ defined from Figure \ref{j13p}.]{\label{j13}\includegraphics[width=0.45\textwidth]{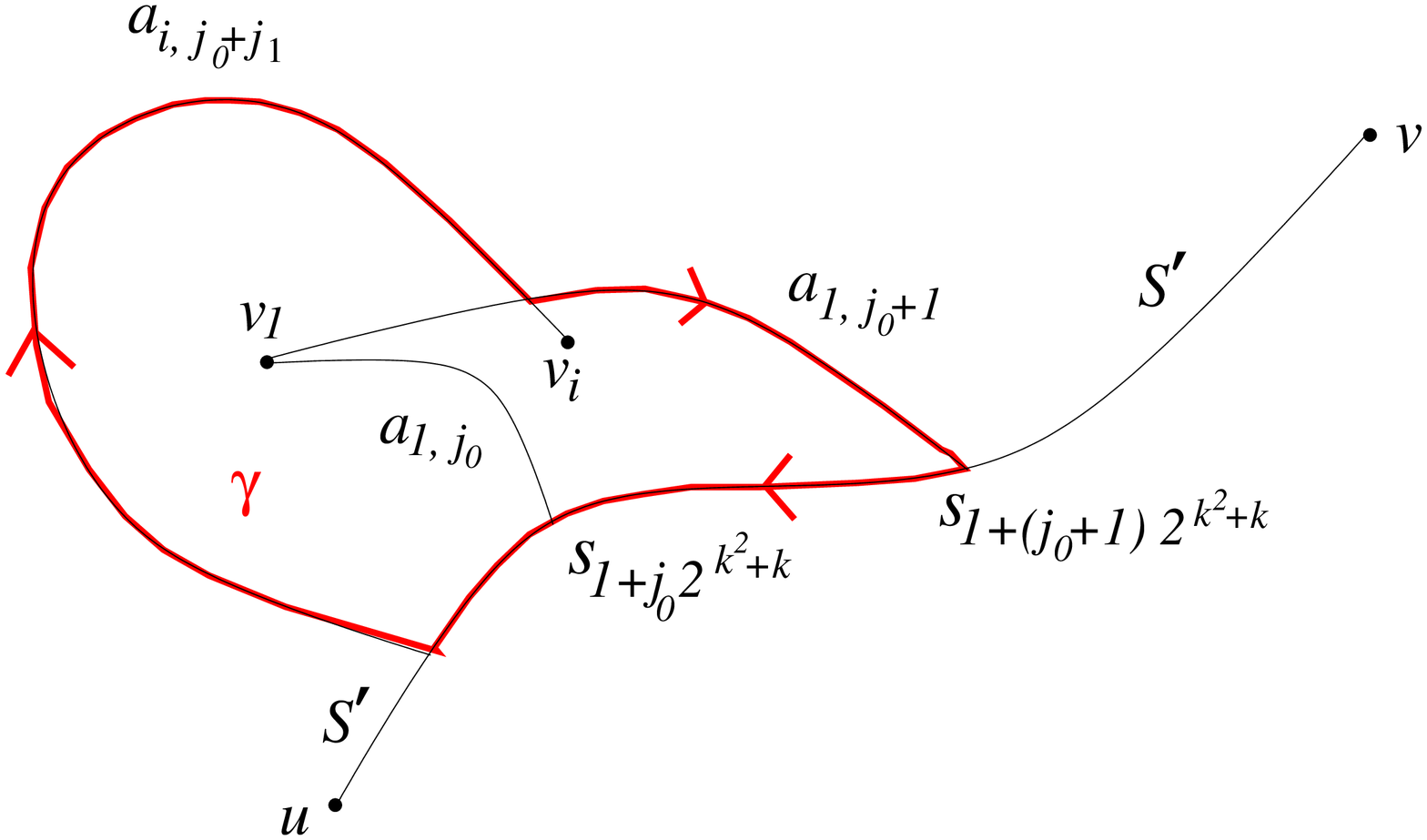}}
                        \caption{Defining $\gamma$ and its orientation.}
\end{figure}

 Since $a_{i,j_0 + j_1}$ intersects $a_{1,j_0 + 1}$, $v_i$ must lie to the right of $\gamma$.  Moreover since the arc from $x$ to $y$ along $a_{1,j_0 + 1}$ is a subset of $\gamma$, the points corresponding to the subsequence

$$ S' = \{s_t \in S \hspace{.2cm}|\hspace{.2cm} 2 + (j_0 + 1)2^{k^2 + k} \leq t \leq (i-1) + (j_0 + j_1)2^{k^2 + k}\}$$

 \noindent lie to the left of $\gamma$.  Hence $\gamma$ separates vertex $v_i$ and the points of $S'$.  Since each arc from $A$ must cross $\gamma$, each arc must cross $a_{1,j_0 + 1}$ since $G$ is simple (these arcs cannot cross $a_{i,j_0 + j_1}$).  See Figure \ref{gamma}.

By the same argument, if the arc $a_{i,j_0 - j_1}$ crosses $a_{1,j_0}$ for $j_1 \geq 1$, then the arcs (emanating out of $v_i$)

$$a_{i,j_0 - 1},a_{i,j_0 - 2},...,a_{i,j_0 - j_1 + 1}$$

\noindent must also cross $a_{1,j_0}$.  Therefore, we have the following observation.

\begin{observation}
\label{obs}
For half of the vertices $v_i \in V'$, the arcs emanating out of $v_i$ satisfy

\begin{enumerate}

\item $a_{i,j_0 + 1}, a_{i,j_0 + 2},...,a_{i,j_0 + 2^k/2}$ all cross $a_{1,j_0 + 1}$, or

\item $a_{i,j_0 - 1}, a_{i,j_0 - 2},...,a_{i,j_0 - 2^k/2}$ all cross $a_{1,j_0}$.

\end{enumerate}

\end{observation}

$\hfill\square$

\noindent Since

$$\frac{|V'|}{2} \geq \frac{2^{k^2 + k} - 1}{2\cdot 2^{k}} \geq 2^{(k-1)^2 + (k-1)},$$

\noindent by Observation \ref{obs} we have an $(2^{(k-1)^2 + (k-1)},2^{k - 1})up$ sequence, whose corresponding arcs all cross either $a_{1,j_0}$ or $a_{1,j_0 +1}$.  By the induction hypothesis, we have $k$ pairwise crossing edges.

$\hfill\square$

Now we are ready to complete the proof of Theorem \ref{simple}.  By Lemma \ref{regular} we know that, say, $S_1$ contains a $2^{k^2 + k}$-regular subsequence of length $|E'|/(4\cdot 2^{k^2 + k})$.  By Theorem \ref{klazar} and Lemma \ref{key}, this subsequence has length at most

$$n2^{k^2 + k}2^{2^{k^2 + 2k} - 3}\left(10\cdot 2^{k^2 + k}\right)^{10\alpha^{2^{k^2 + 2k}}(n)}.$$

\noindent Therefore

$$\frac{2m}{10^4\cdot 4\cdot 2^{k^2 + k}\log^2n} \leq \frac{|E'|}{4\cdot 2^{k^2 + k}} \leq n2^{k^2 + k}2^{2^{k^2 + 2k} - 3}\left(10\cdot 2^{k^2 + k}\right)^{10\alpha^{2^{k^2 + 2k}}(n)}$$

\noindent which implies

$$m \leq 4\cdot 10^4\cdot 2^{2k^2 + 2k}2^{2^{k^2 + 2k} - 3} n\left(10\cdot 2^{k^2 + k}\right)^{10\alpha^{2^{k^2 + 2k}}(n)}\log^2n.$$

\noindent Since $c_k = 10^5\cdot 2^{k^2 + 2k}$, $\alpha(n) \geq 2$ and $k \geq5$, we have

$$m \leq (n\log^2 n)2^{\alpha^{c_k}(n)}.$$

$\hfill\square$

\section{$x$-monotone}

In this section we will prove Theorem \ref{xmono}.

\medskip

\noindent \textbf{Proof of Theorem \ref{xmono}.}   For $k \geq 2$, let $g_k(n)$ be the maximum number of edges in a (not simple) $k$-quasi-planar graph whose edges are drawn as $x$-monotone curves.  We will prove by induction on $n$ that

  $$g_k(n) \leq 2^{c k^3}n\log n$$

  \noindent where $c$ is a sufficiently large absolute constant.  The base case is trivial.  For the inductive step, let $G = (V,E)$ be a $k$-quasi-planar topological graph whose edges are drawn as $x$-monotone curves, and let the vertices be labeled $1,2,...,n$.  Then let $L$ be the vertical line that partitions the vertices into two parts, $V_1$ and $V_2$, such that $|V_1| = \lfloor n/2\rfloor$ vertices lie to the left of $L$, and $|V_2| = \lceil n/2\rceil$ vertices lie to the right of $L$.  Furthermore, let $E_1$ denote the set of edges induced by $V_1$, $E_2$ be the set of edges induced by $V_2$, and $E'$ be the set of edges that intersect $L$. Clearly, we have

$$|E_1| \leq g_k(\lfloor n/2\rfloor) \hspace{1cm}\textnormal{and}\hspace{1cm} |E_2| \leq g_k(\lceil n/2\rceil).$$

\noindent Hence it suffices that show that

\begin{equation}
\label{goal}
|E'|\leq 2^{c k^3/2}n,
\end{equation}

\noindent since this would imply

$$g_k(n) \leq g_k(\lfloor n/2\rfloor) + g_k(\lceil n/2\rceil)+ 2^{c k^3/2}n \leq  2^{c k^{3}} n\log n. $$

 For the rest of the proof, we will only consider the edges from $E'$.  Now for each vertex $v_i \in V_1$, consider the graph $G_i$ whose vertices are the edges with $v_i$ as a left endpoint, and two vertices in $G_i$ are adjacent if the corresponding edges cross at some point to the left of $L$.  Since $G_i$ is an \emph{incomparability graph} (see \cite{dilworth}, \cite{fox2}) and does not contain a clique of size $k$, $G_i$ contains an independent set of size $|E(G_i)|/(k-1)$.  We keep all edges that correspond to the elements of this independent set, and discard all other edges incident to $v_i$.  After repeating this process on all vertices in $V_1$, we are left with at least $|E'|/(k-1)$ edges.

Now we continue this process on the other side.  For each vertex $v_j \in V_2$, consider the graph $G_j$ whose vertices are the edges with $v_j$ as a right endpoint, and two vertices in $G_j$ are adjacent if the corresponding edges cross at some point to the right of $L$.  Since $G_j$ is an incomparability graph and does not contain a clique of size $k$, $G_j$ contains an independent set of size $|E(G_j)|/(k-1)$.  We keep all edges that corresponds to this independent set, and discard all other edges incident to $v_j$.  After repeating this process on all vertices in $V_2$, we are left with at least $|E'|/(k-1)^2$ edges.

We order the remaining edges $e_1,e_2,...,e_m$ in the order in which they intersect $L$ from bottom to top.  We define two sequences $S_1 = a_1,a_2,...,a_m$ and $S_2 = b_1,b_2,...,b_m$ such that $a_i$ denotes the left endpoint of edge $e_i$ and $b_i$ denotes the right endpoint of $e_i$.  Now we need the following lemma.

\begin{lemma}
\label{key2}
Neither of the sequences $S_1$ or $S_2$ contains a subsequence of type up-down-up$(k^3+2)$.

\end{lemma}

\noindent \textbf{Proof.}  By symmetry, it suffices to show that $S_1$ does not contain a subsequence of type \emph{up-down-up}$(k^3 + 2)$.  For the sake of contradiction, suppose $S_1$ did contain a subsequence of type \emph{up-down-up}$(k^3 + 2)$.  Then there is a sequence

$$S = s_1,s_2,...,s_{3(k^3 + 2)-2}$$

\noindent such that the integers $s_1,...,s_{k^3 + 2}$ are pairwise different and for $i = 1,2,...,k^3 + 2$ we have

$$s_i = s_{2(k^3 + 2) - i} = s_{2(k^3 + 2) - 2 + i}.$$

For each $i = 1,2,...,k^3 + 2$, let $v_i \in V_1$ denote the label (vertex) of $s_i$ and let $x_i$ denote the $x$-coordinate of vertex $v_i$.  Moreover, let $a_i$ be the arc emanating from vertex $v_i$ to the point on $L$ that corresponds to $s_{2(k^3+2) - i}$.  Note that the set of arcs $A = \{a_2,a_3,...,a_{k^3 + 1}\}$ are ordered downwards along $L$, and corresponds to the ``middle" part of the up-down-up sequence.  We define two partial orders on $A$ as follows.

$$\begin{array}{cccccc}
   a_i \prec_1 a_j & \textnormal{if} &  i < j, &x_i < x_j & \textnormal{and the arcs $a_i,a_j$ do not intersect,}\\\\
   a_i \prec_2 a_j & \textnormal{if} &  i < j, &x_i > x_j & \textnormal{and the arcs $a_i,a_j$  do not intersect.}
 \end{array}$$

Clearly, $\prec_1$ and $\prec_2$ are partial orders.  If two arcs are not comparable by either $\prec_1$ or $\prec_2$, then they must cross.  Since $G$ does not contain $k$ pairwise crossing edges, by Dilworth's Theorem, there exist $k$ arcs $\{a_{i_1},a_{i_2},...,a_{i_k}\}$ such that they are pairwise comparable by either $\prec_1$ or $\prec_2$.  Now the proof falls into two cases.

\medskip

\begin{figure}[h]
  \centering
\subfigure{ \includegraphics[width=0.27 \textwidth]{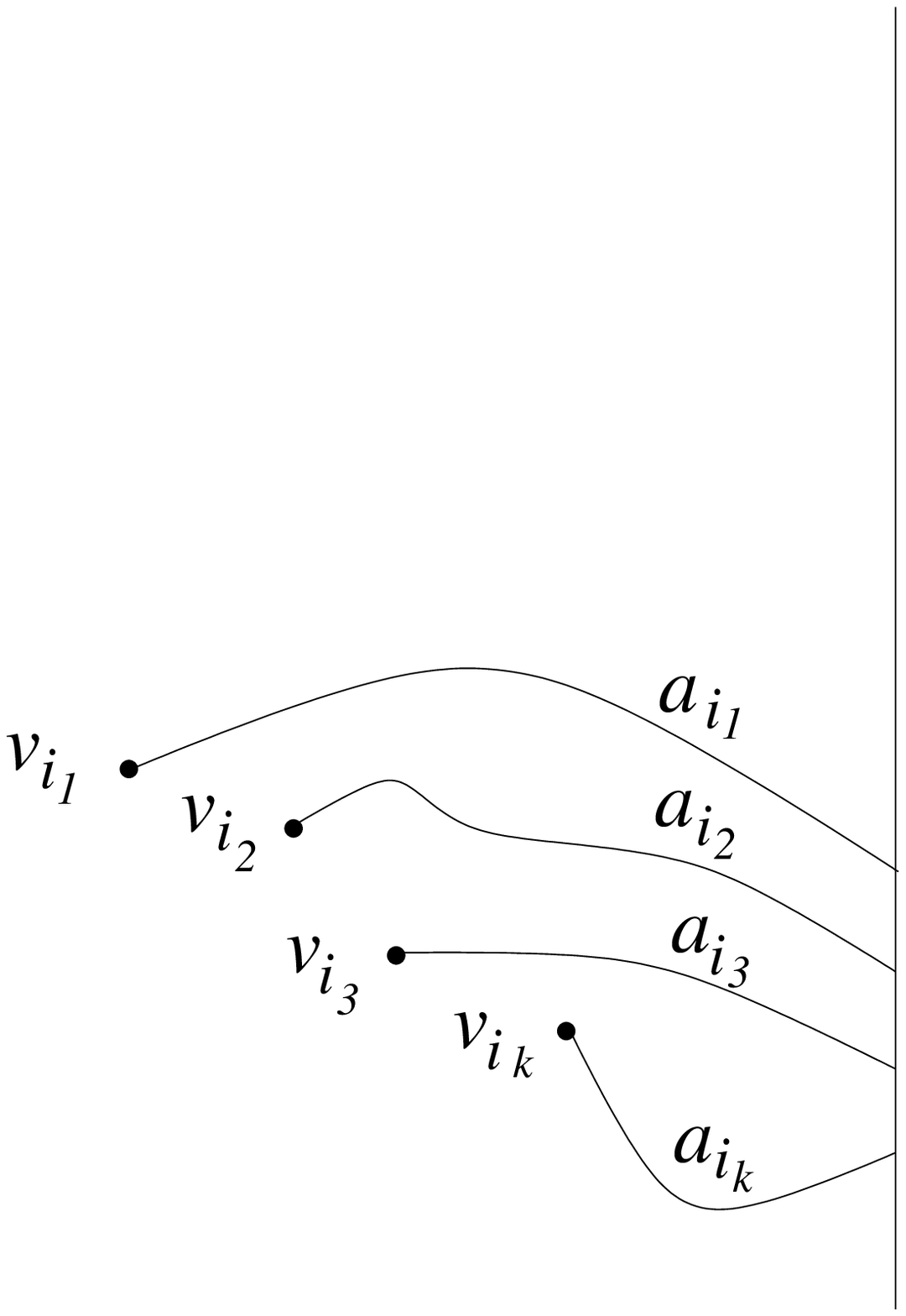}} \hspace{2cm}
\subfigure{ \includegraphics[width=0.35\textwidth]{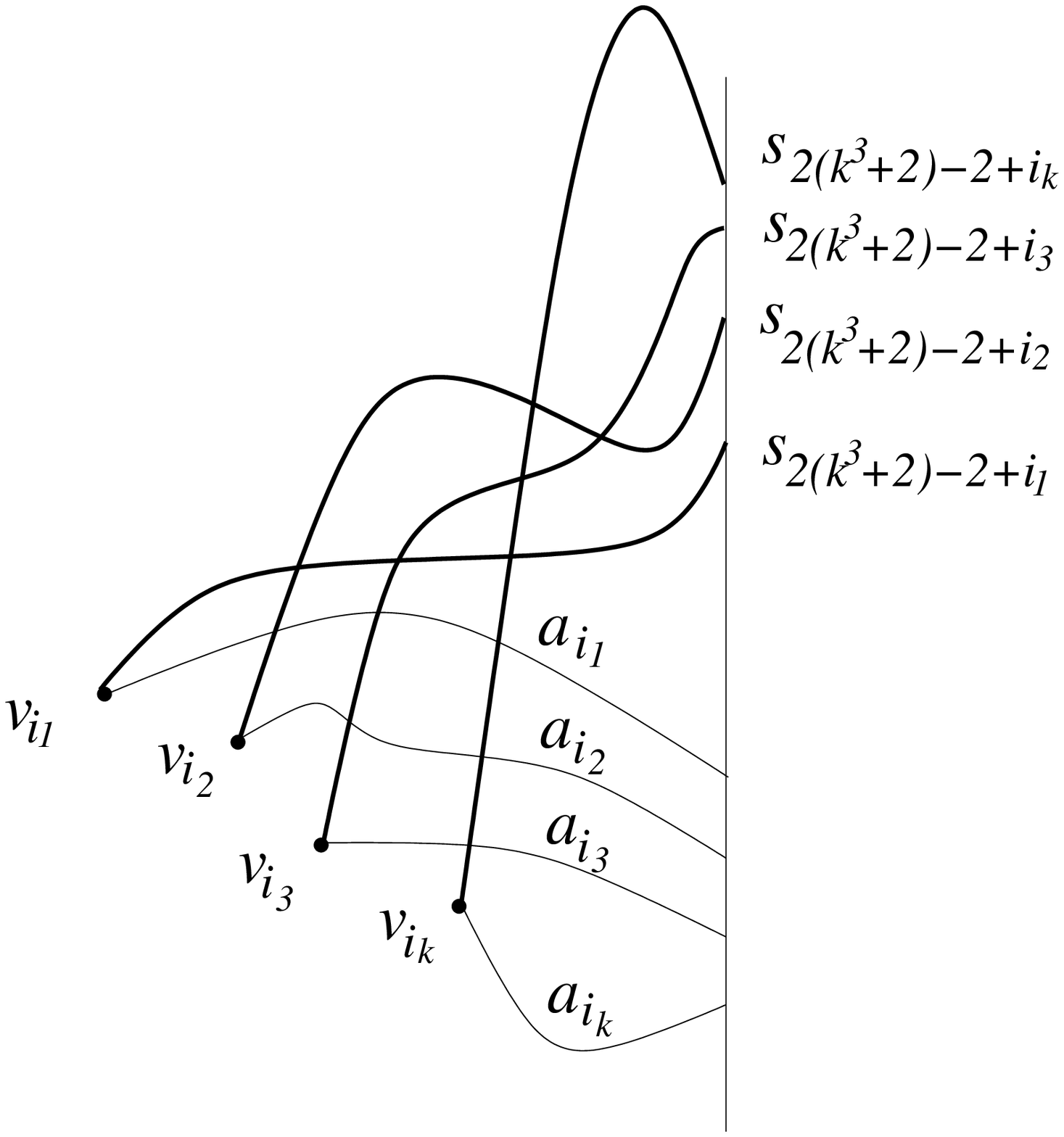}}
                        \caption{Case 1.}
  \label{prec1}
\end{figure}

\noindent \emph{Case 1.}  Suppose that $a_{i_1} \prec_1 a_{i_2} \prec_1 \cdots \prec_1 a_{i_k}$.  Then the arcs emanating from $v_{i_1},v_{i_2},...,v_{i_k}$ to the points corresponding to $s_{2(k^3+2)-2 + i_1},s_{2(k^3 + 2)-2 + i_2},...,s_{2(k^3 + 2)-2 + i_k}$ are pairwise crossing.  See Figure \ref{prec1}.

\medskip

\begin{figure}[h]
  \centering
\subfigure{ \includegraphics[width=0.28\textwidth]{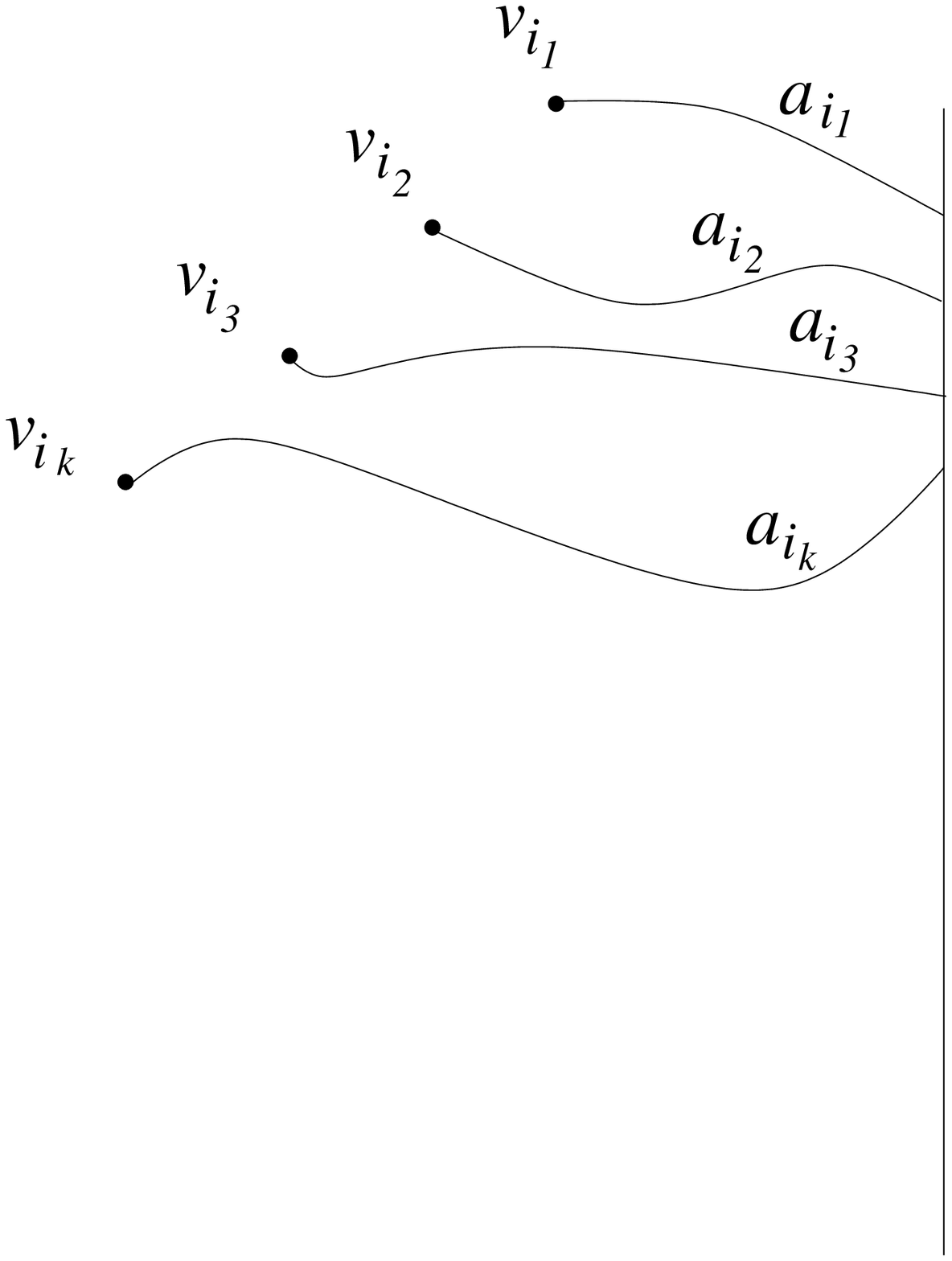}} \hspace{2cm}
\subfigure{ \includegraphics[width=0.35\textwidth]{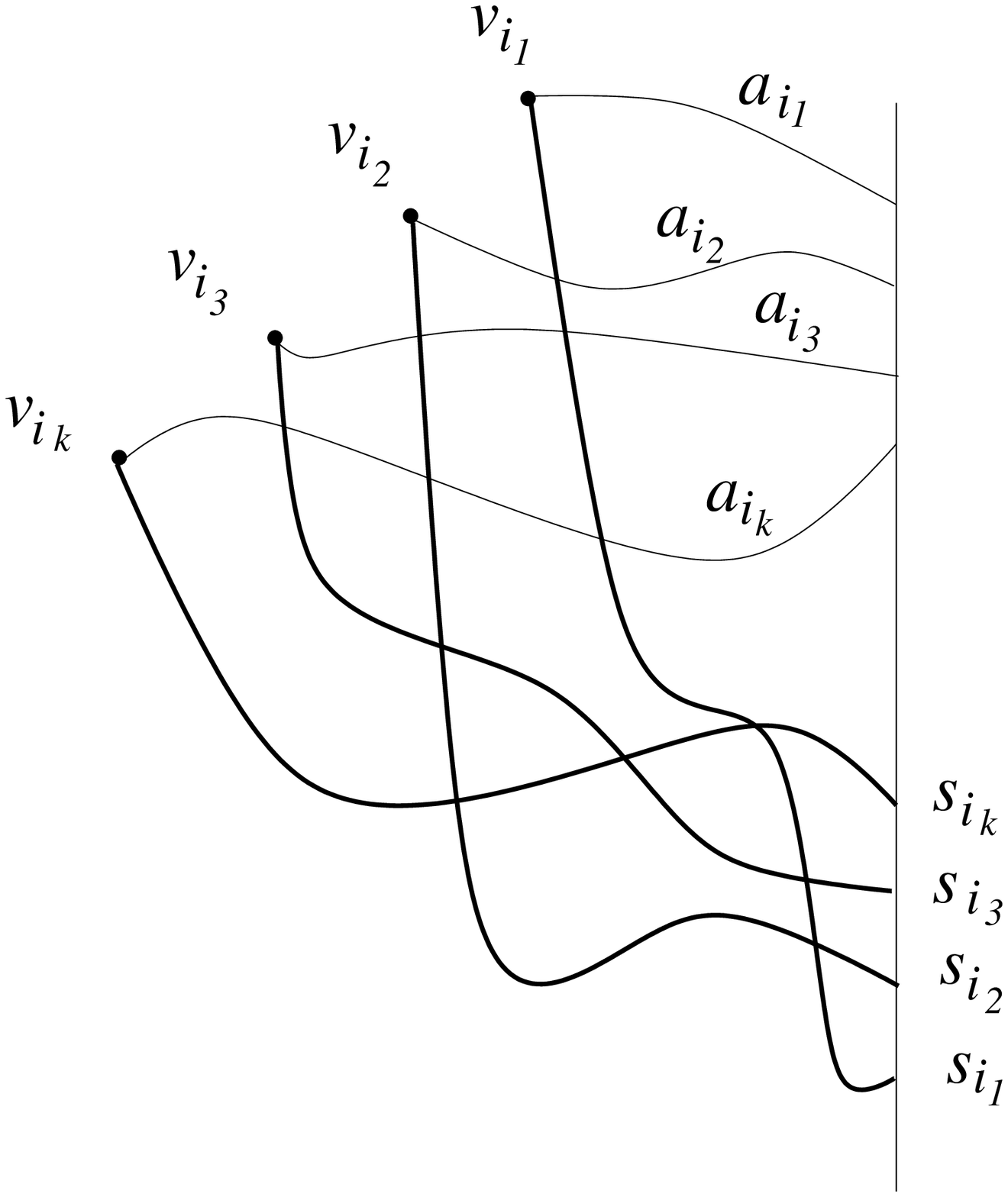}}
                        \caption{Case 2.}
  \label{prec2}
\end{figure}

\noindent \emph{Case 2.}    Suppose that $a_{i_1} \prec_2 a_{i_2} \prec_2 \cdots \prec_2 a_{i_k}$. Then the arcs emanating from $v_{i_1},v_{i_2},...,v_{i_k}$ to the points corresponding to $s_{i_1},s_{i_2},...,s_{i_k}$ are pairwise crossing.  See Figure \ref{prec2}.

$\hfill\square$

We are now ready to complete the proof of Theorem \ref{xmono}.  By Lemma \ref{regular}, we know that, say, $S_1$ contains a $(k^3+2)$-regular subsequence of length

$$\frac{|E'|}{4(k^3 + 2)(k-1)^2}.$$

\noindent By lemma \ref{updownup} and \ref{key2}, this subsequence has length at most $2^{c' k^3}n$, where $c'$ is an absolute constant.  Hence

$$\frac{|E'|}{4(k^3 + 2)(k-1)^2} \leq 2^{c' k^3}n$$

\noindent implies

$$|E'| \leq 4k^52^{c'k^3}n \leq 2^{c k^3/2} n$$

\noindent for a sufficiently large absolute constant $c$.

$\hfill\square$

\end{document}